\documentclass[a4paper,11pt]{article}
\usepackage{latexsym}
\usepackage{amssymb}
\usepackage{amsfonts}
\usepackage{amsmath}
\usepackage{indentfirst}
\usepackage{graphicx}
\usepackage{epsfig}
\usepackage{slashbox,colortbl}
\usepackage[enableskew]{youngtab}
\usepackage{fullpage}

\newtheorem{prop}{Proposition}[section]

\newcommand{\kinN}{k\in \mathbf{N}}

\newcommand{\cvd}{\hfill $\blacksquare$\bigskip}
\date{}

\author{Antonio Bernini\thanks{Dipartimento di Matematica e Informatica ``U. Dini", viale Morgagni 65, University
of Firenze, Firenze, Italy, \tt{antonio.bernini@unifi.it,
luca.ferrari@unifi.it, pinzani@dsi.unifi.it}}, Luca
Ferrari$^\dag$, Renzo Pinzani$^\dag$, Julian West\thanks{Heilbronn
Institute for Mathematical Research, University of Bristol,
Bristol, UK, \tt{julian.ndp@gmail.com}}}

\title{The Dyck pattern poset\footnote{LF and RP partially supported by INdAM project: \emph{Propriet\'a algebriche e
combinatorie dei cammini discreti}}}

\begin{document}

\maketitle

\begin{abstract}
We introduce the notion of \emph{pattern} in the context of
lattice paths, and investigate it in the specific case of Dyck
paths. Similarly to the case of permutations, the
pattern-containment relation defines a poset structure on the set
of all Dyck paths, which we call the \emph{Dyck pattern poset}.
Given a Dyck path $P$, we determine a formula for the number of
Dyck paths covered by $P$, as well as for the number of Dyck paths
covering $P$. We then address some typical pattern-avoidance
issues, enumerating some classes of pattern-avoiding Dyck paths.
Finally, we offer a conjecture concerning the asymptotic behavior
of the sequence counting Dyck paths avoiding a generic pattern and
we pose a series of open problems regarding the structure of the
Dyck pattern poset.
\end{abstract}

\section{Introduction}

One of the most investigated and fruitful notions in contemporary
combinatorics is that of a \emph{pattern}. Historically it was
first considered for permutations \cite{Kn}, then analogous
definitions were provided in the context of many other structures,
such as set partitions \cite{Go,Kl,Sa}, words \cite{Bj,Bu}, and
trees \cite{DPTW,Gi,R}. Perhaps all of these examples have been
motivated or informed by the more classical notion of graphs and
subgraphs. Informally speaking, given a specific class of
combinatorial objects, a pattern can be thought of as an
occurrence of a small object inside a larger one; the word
``inside" means that the pattern is suitably embedded into the
larger object, depending on the specific combinatorial class of
objects. The main aim of the present work is to introduce the
notion of pattern in the context of lattice paths and to begin its
systematic study in the special case of Dyck paths.

\bigskip

For our purposes, a \emph{lattice path} is a path in the discrete
plane starting at the origin of a fixed Cartesian coordinate
system, ending somewhere on the $x$-axis, never going below the
$x$-axis and using only a prescribed set of steps $\Gamma$. We
will refer to such paths as \emph{$\Gamma$-paths}. This definition
is extremely restrictive if compared to what is called a lattice
path in the literature, but it will be enough for our purposes.
Observe that a $\Gamma$-path can be alternatively described as a
finite word on the alphabet $\Gamma$ obeying certain conditions.
Using this language, we say that the \emph{length} of a
$\Gamma$-path is simply the length of the word which encodes such
a path. Among the classical classes of lattice paths, the most common are those
using only steps $U(p)=(1,1)$, $D(own)=(1,-1)$ and $H(orizontal)=(1,0)$;
with these definitions, Dyck, Motzkin and Schr\"oder paths correspond
respectively to the set of steps $\{ U,D\}$, $\{ U,H,D\}$ and $\{
U,H^2 ,D\}$.

\bigskip

Consider the class $\mathcal{P}_\Gamma$ of all $\Gamma$-paths, for
some choice of the set of steps $\Gamma$. Given $P,Q\in
\mathcal{P}_\Gamma$ having length $k$ and $n$ respectively, we say
that \emph{$Q$ contains (an occurrence of) the pattern $P$}
whenever $P$ occurs as a subword of $Q$. So, for instance, in the
class of Dyck paths, $UUDUDDUDUUDD$ contains the pattern $UUDDUD$,
whereas in the class of Motzkin paths, $UUHDUUDHDDUDHUD$ contains
the pattern $UHUDDHUD$. When $Q$ does not contain any occurrence
of $P$ we will say that $Q$ \emph{avoids} $P$. In the Dyck case,
the previously considered path $UUDUDDUDUUDD$ avoids the pattern
$UUUUDDDD$.

\begin{figure}[h!]
\begin{center}
    \includegraphics[scale=.4]{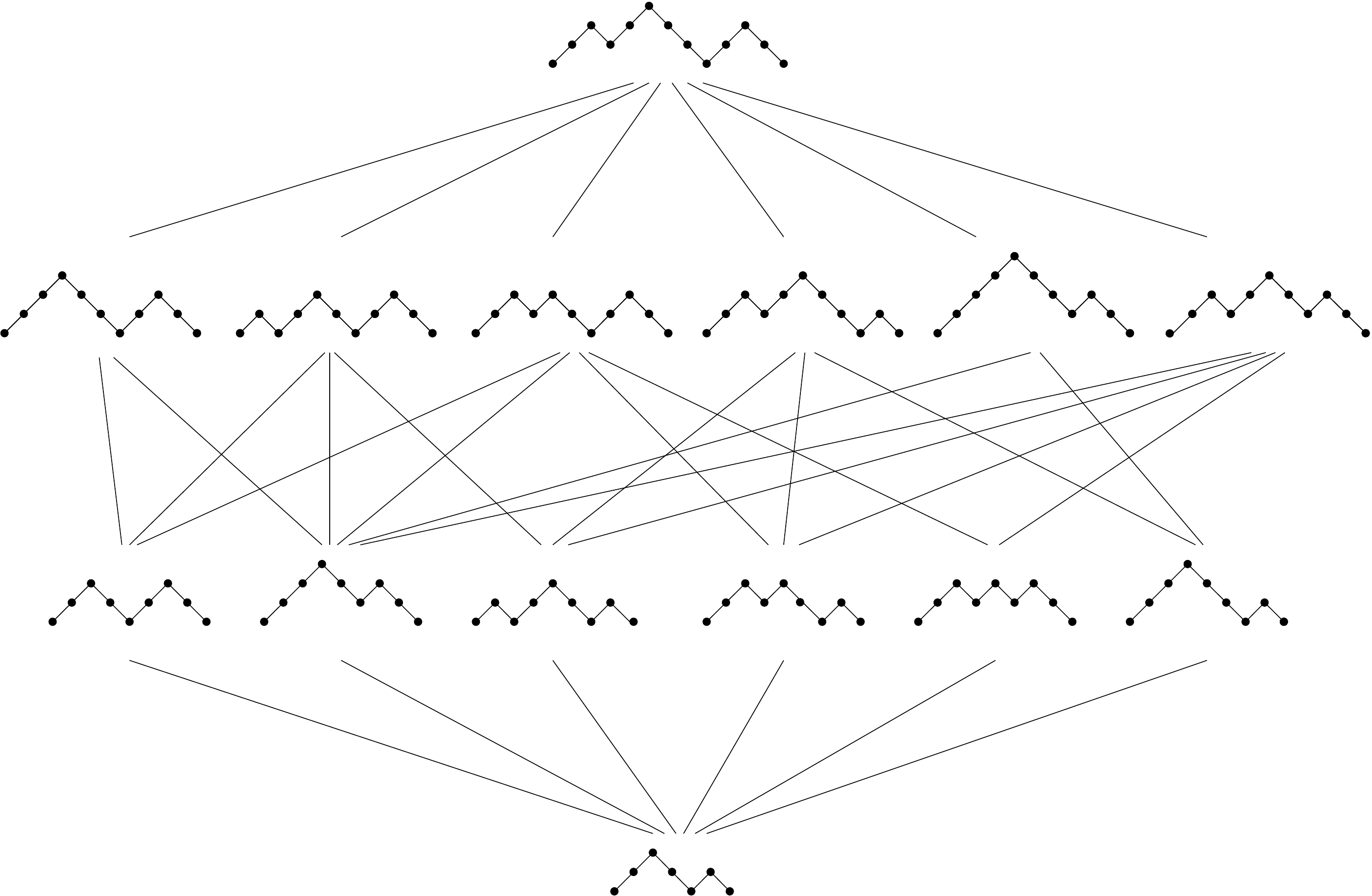}
\end{center}
\caption{An interval of rank 3 in the Dyck pattern poset.}
\label{interval}
\end{figure}

This notion of pattern gives rise to a partial order in a very
natural way, by declaring $P\leq Q$ when $P$ occurs as a pattern
in $Q$. In the case of Dyck paths, the resulting poset will be
denoted by $\mathcal{D}$. It is immediate to notice that
$\mathcal{D}$ has a minimum (the empty path), does not have a
maximum, is locally finite and is ranked (the rank of a Dyck path
is given by its semilength). As an example, in Figure
\ref{interval} we provide the Hasse diagram of an interval in the
Dyck pattern poset.

Observe that this notion of pattern for paths is very close to the
analogous notion for words (considered, for instance, in
\cite{Bj}, where the author determines the M\"obius function of
the associated pattern poset). Formally, instead of considering
the set of \emph{all} words of the alphabet $\{ U,D\}$, we
restrict ourselves to the set of Dyck words (so what we actually
do is to consider a subposet of Bj\"orner's poset). However, the
conditions a word has to obey in order to belong to this subposet
(which translate into the fact of being a Dyck word) make this
subposet highly nontrivial, and fully justify our approach,
consisting of the study of its properties independently of its
relationship with the full word pattern poset.

\section{The covering relation in the Dyck pattern poset}

\begin{figure}[h]

\begin{center}
    \includegraphics[scale=.5]{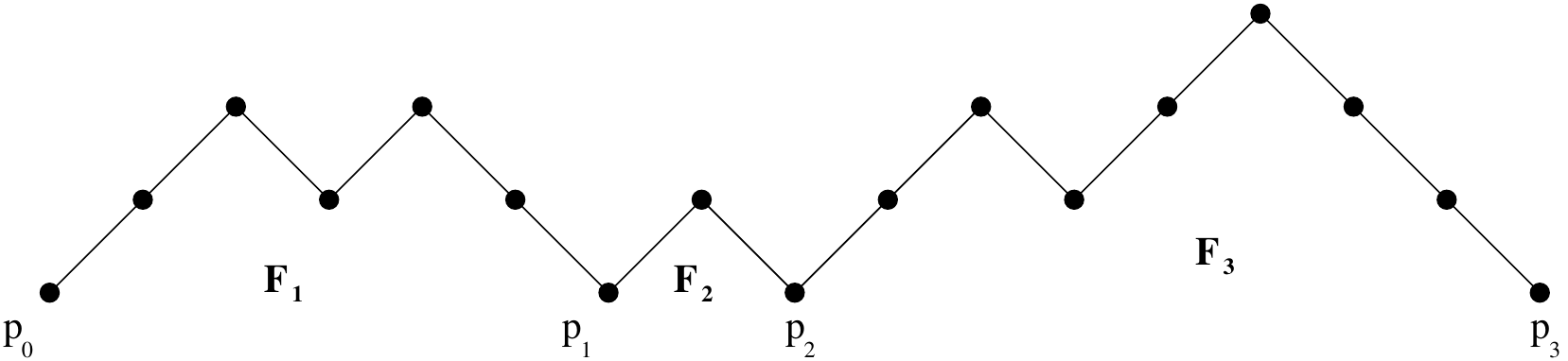}
\end{center}
\caption{A Dyck path having three factors.} \label{notation}
\end{figure}

In the Dyck pattern poset $\mathcal{D}$, following the usual
notation for covering relation, we write $P\prec Q$ ($Q$
\emph{covers} $P$) to indicate that $P\leq Q$ and the rank of $P$
is one less than the rank of $Q$ (i.e., $rank(P)=rank(Q)-1$). Our
first result concerns the enumeration of Dyck paths covered by a
given Dyck path $Q$. We need some notation before stating it. Let
$k+1$ be the number of points of $Q$ lying on the $x$-axis (call
such points $p_0,p_1,\ldots,p_k$). Then $Q$ can be factorized into
$k$ Dyck factors $F_1,\ldots,F_k$, each $F_i$ starting at
$p_{i-1}$ and ending at $p_i$. Let $n_i$ be the number of ascents
in $F_i$ (an ascent being a consecutive run of $U$ steps; $n_i$
also counts both the number of descents and the number of peaks in
$F_i$). Moreover, we denote by $|UDU|$ and $|DUD|$ the number of
occurrences in a Dyck path of a consecutive factor $UDU$ and
$DUD$, respectively. In the path $Q$ of Figure \ref{notation}, we
have $n_1=2$, $n_2=1$, $n_3=2$, $|UDU|=3$, and $|DUD|=2$.

\begin{prop}
If $Q$ is a Dyck path with $k$ factors $F_1 ,\ldots F_k$, with
$F_i$ having $n_i$ ascents, then the number of Dyck paths covered
by $Q$ is given by
\begin{equation}\label{covered}
\frac{\sum_{i=1}^k{n_i}^2+(\sum_{i=1}^k
n_i)^2}{2}-|UDU|-|DUD|\quad .
\end{equation}
\end{prop}

\emph{Proof.} We proceed by induction on $k$. If $Q$ is any Dyck
path having only one factor (and so necessarily $n_1$ ascents),
then a path $P$ such that $P\prec Q$ is obtained by choosing (and
then removing) a $U$ step and a $D$ step from an ascent and a
descent of $Q$, respectively. This can be done in $n_1^2$
different ways. Note that, once an ascent and a descent have been
selected, the path $P$ is uniquely determined if each of the two
steps $U$ and $D$ is chosen within the same ascent and descent,
respectively. Moreover, for each $UDU$ (resp., $DUD$) occurring in
$Q$, removing the $D$ step from the $UDU$ (resp., the $U$ step
from the $DUD$) and a $U$ (resp., $D$) step from the ascent
(resp., descent) either immediately before $D$ (resp., $U$) or
immediately after $D$ (resp., $U$) produces the same path $P$
covered by $Q$. Therefore, these paths would be counted twice if
the term $n_1^2$ were not corrected by subtracting both $|UDU|$
and $|DUD|$. This leads to formula (\ref{covered}) in the case
$k=1$.

Now suppose that $\tilde{Q}$ is a Dyck path which has $k>1$
factors $F_1,\ldots,F_k$, each factor $F_i$ having $n_i$ ascents.
Let $l$ be the total number of $UDU$ and $DUD$ (i.e.
$l=|UDU|+|DUD|$) in $\tilde{Q}$. If a new factor $F_{k+1}$ having
$n_{k+1}$ ascents and a total number $l_{k+1}$ of $UDU$ and $DUD$
factors is appended to $\tilde{Q}$ (after $F_k$), then the paths
covered by the new path $Q$ can be obtained by removing a $D$ step
and a $U$ step either both belonging to $\tilde{Q}$, or both
belonging to $F_{k+1}$, or one belonging to $\tilde{Q}$ and the
other one belonging to $F_{k+1}$.

We start by supposing that the two factors $F_k$ and $F_{k+1}$ are
both different from $UD$. In the first of the above cases, the
number of covered paths is given by formula (\ref{covered}) thanks
to our inductive hypothesis (since the removal of the steps $U$
and $D$ involves only the first $k$ factors of the Dyck path). The
second case is easily dealt with using the induction hypothesis as
well, namely applying the base case ($k=1$) to the last factor
$F_{k+1}$. Finally, concerning the last case, notice that the step
$D$ must be removed from $\tilde{Q}$, and the step $U$ must be
removed from $F_{k+1}$, otherwise the resulting path would fall
below the $x$-axis. Then, the $D$ step can be selected from
$\sum_{i=1}^kn_i$ different descents of $\tilde{Q}$, while the $U$
step can be chosen among the steps of the $n_{k+1}$ ascents of
$F_{k+1}$, leading to $n_{k+1}\cdot \sum_{i=1}^{k}n_i$ different
paths covered by $Q$. Summing the contributions of the three
cases considered above, we obtain:

\begin{eqnarray}\label{partial}
& &\frac{\sum_{i=1}^k{n_i}^2+(\sum_{i=1}^k n_i)^2}{2}-l+n_{k+1}^2-l_{k+1}+n_{k+1}\sum_{i=1}^kn_i \nonumber\\
&=&\frac{\sum_{i=1}^{k+1}{n_i}^2+(\sum_{i=1}^{k+1}
n_i)^2}{2}-l-l_{k+1}\quad.
\end{eqnarray}

However, we still have to take into account the cases in which
$F_k$ and/or $F_{k+1}$ are equal to $UD$. If $F_k =F_{k+1}=UD$,
then in formula (\ref{partial}) we have to subtract $2$ (since we
have one more factor $UDU$ and one more factor $DUD$ than those
previously counted). In the remaining cases, there is only one
more factor (either $UDU$ or $DUD$), thus in formula
(\ref{partial}) we have to subtract $1$. In all cases, what we get
is precisely formula (\ref{covered}).\cvd

In a similar fashion, we are also able to find a formula for the
number of all Dyck paths which cover a given path.

\begin{prop} If $Q$ is a Dyck path of semilength $n$ with $k$ factors $F_1 ,\ldots F_k$, with
$F_i$ having semilength $f_i$, then the number of Dyck paths
covering $Q$ is given by
\begin{equation}\label{covering}
1+\sum_{i}f_i ^2 +\sum_{i<j}f_i f_j .
\end{equation}
\end{prop}

\emph{Proof.}\quad A path $P$ covers $Q$ if and only if it is
obtained from $Q$ by suitably inserting an up step $U$ and a down
step $D$. Thus the set of all Dyck paths covering $Q$ can be
determined by choosing, in all possible ways, two positions
(inside $Q$) in which to insert an up step and a down step.
Clearly, in performing these insertions, we must take care not to
fall below the $x$-axis.

Let $P$ cover $Q$ and denote with $R$ the first occurrence (from
the left) of $Q$ in $P$. There are precisely two steps
in $P$ (a $U$ and a $D$) which do not belong to $R$. We
distinguish three distinct cases.

\begin{enumerate}

\item The last step of $R$ is the third-to-last step of $P$ (so that
$R$ is a prefix of $P$). This means that the two added steps
are the last two steps of $P$ (which therefore ends with $UD$),
and it is clear that there is precisely one path $P$ covering $Q$
which falls into this case.

\item The last step of $R$ is the second-to-last step of $P$. This
means that the $D$ step inserted into $R$ is the last step of $P$.
Thus $P$ is obtained by inserting a $U$ step somewhere in $R$
(except that at the end, since in this case we will return to the
previous case). The number of path $P$ of this form is then given
by the number of different places of $R$ in which we are allowed
to insert a new up step. Since $R$ is required to be the first
occurrence of $Q$ in $P$, it can be shown that a new up step can
be inserted immediately before each down step of $R$. There are
precisely $n$ Dyck paths of this form.

\item The last step of $R$ is the last step of $P$. In this case,
$P$ is obtained from $R$ by suitably inserting an up step and a
down step. We can consider two distinct cases. If $U$ and $D$ are
inserted into the same factor of $R$, then we can either insert
$U$ before $D$ or vice versa. In this specific case, since $R$ has
to be the first occurrence of $Q$ inside $P$, $U$ can be inserted
immediately before each $D$ step of the factor and $D$ can be
inserted immediately before each $U$; similarly, $D$ can be
inserted immediately before each $U$, except of course for the
very first step of the factor, moreover $D$ can also be inserted
at the end of the factor. There is however one factor that behaves
in a slightly different way. If we choose to insert the two new
steps into the last factor of $P$, then we cannot insert a $D$ at
the end of the factor (since we are supposing that the last step
of $R$ is also the last step of $P$). Thus, if we insert $U$ and
$D$ into the factor $F_i$, $i<k$, then we obtain $f_i ^2$
different paths $P$ of this form, whereas if we insert $U$ and $D$
into $F_k$ we get a total of $f_k (f_k -1)$ paths. So, in this
specific case, the total number of paths thus obtained is
$\sum_{i=1}^{k}f_i ^2 -f_k$. On the other hand, if we choose to
insert $U$ and $D$ into two distinct factors, then $U$ must be
inserted before $D$ (otherwise the resulting path would fall below
the $x$-axis). If we decide to insert $D$ into the factor $F_i$,
$i<k$ (for which, by an argument similar to the above one, we have
$f_i$ possibilities), then we can insert $U$ into any of the
preceding factors, whence in $\sum_{j=1}^{i-1}f_j$ ways. Instead,
if $D$ is inserted into $F_k$, we only have $f_k -1$
possibilities, and we can then insert $U$ in any of the first
$k-1$ factors, for a total of $\sum_{j=1}^{k-1}f_j$ different
paths thus obtained. Thus, in this last case, the total number of
paths $P$ having this form is given by $\sum_{i=1}^{k-1}\left( f_i
\cdot \sum_{j=1}^{i-1}f_j \right) +(f_k -1)\cdot
\sum_{j=1}^{k-1}f_j$.

\end{enumerate}

Finally, summing up all the quantities obtained so far, we find
the following expression for the number of paths covering a given
path $Q$:

\begin{eqnarray*}
& &1+n+\sum_{i=1}^{k}f_i ^2 -f_k +\sum_{i=1}^{k}\left( f_i \cdot
\sum_{j=1}^{i-1}f_j \right) -\sum_{i=1}^{k-1}f_i \nonumber
\\ &=&1+\sum_{i=1}^{k}f_i ^2 +\sum_{i<j}f_i f_j .
\end{eqnarray*}

This is precisely formula (\ref{covering}).\cvd

\section{Enumerative results on pattern avoiding Dyck paths}

In the present section we will be concerned with the enumeration
of some classes of pattern avoiding Dyck paths. Similarly to what
has been done for other combinatorial structures, we are going to
consider classes of Dyck paths avoiding a single pattern, and we
will examine the cases of short patterns. Specifically, we will
count Dyck paths avoiding any single path of length $\leq 3$; each
case will arise as a special case of a more general result
concerning a certain class of patterns.

Given a pattern $P$, we denote by $D_n (P)$ the set of all Dyck
paths of semilength $n$ avoiding the pattern $P$, and by $d_n (P)$
the cardinality of $D_n (P)$.

\subsection{The pattern $(UD)^k$}

This is one of the easiest cases.

\begin{prop} For any $\kinN$, $Q\in D_n ((UD)^k)$ if and only if
$Q$ has at most $k-1$ peaks.
\end{prop}

\emph{Proof.}\quad A Dyck path
$Q=U^{a_1}D^{b_1}U^{a_2}D^{b_2}\cdots U^{a_h}D^{b_h}$ contains the
pattern $(UD)^k$ if and only if $h\geq k$, that is $Q$ has at
least $k$ peaks.\cvd

Since it is well known that the number of Dyck paths of semilength
$n$ and having $k$ peaks is given by the Narayana number $N_{n,k}$
(sequence A001263 in \cite{Sl}), we have that $d_n
((UD)^k)=\sum_{i=0}^{k-1}N_{n,i}$ (partial sums of Narayana
numbers). Thus, in particular:

\begin{itemize}

\item[-] $d_n (UD)=0$;

\item[-] $d_n (UDUD)=1$;

\item[-] $d_n (UDUDUD)=1+{n\choose 2}$.

\end{itemize}

\subsection{The pattern $U^{k-1}DUD^{k-1}$}

Let $Q$ be a Dyck path of length $2n$ and $P=U^{k-1}DUD^{k-1}$.
Clearly if $n<k$, then $Q$ avoids $P$, and if $n=k$, then all Dyck
paths of length $2n$ except one ($Q$ itself) avoid $Q$. Therefore:

\begin{itemize}

\item $d_n(P)=C_n$ if $n<k$, and

\item $d_n(P)=C_n -1$ if $n=k$,

\end{itemize}

\noindent where $C_n$ is the $n$-th Catalan number.

Now suppose $n>k$. Denote by $A$ the end point of the $(k-1)$-th
$U$ step of $Q$. It is easy to verify that $A$ belongs to the line
$r$ having equation $y=-x+2k-2$. Denote with $B$ the starting
point of the $(k-1)$-th-to-last $D$ step of $Q$. An analogous
computation shows that $B$ belongs to the line $s$ having equation
$y=x-\left(2n-2k+2\right)$.

Depending on how the two lines $r$ and $s$ intersect, it is
convenient to distinguish two cases.

\begin{enumerate}

\item If $2n-2k+2\geq2k-4$ (i.e. $n\geq 2k-3$), then $r$ and $s$
intersect at height $\leq 1$, whence $x_A\leq x_B$ (where $x_A$
and $x_B$ denote the abscissas of $A$ and $B$, respectively). The
path $Q$ can be split into three parts (see Figure
\ref{av_UUDUDD_}): a prefix $Q_A$ from the origin $(0,0)$ to
$A$, a path $X$ from $A$ to $B$, and a suffix $Q_B$ from $B$ to
the last point $(2n,0)$.

\begin{figure}[h!]
\begin{center}
\includegraphics[scale=.3]{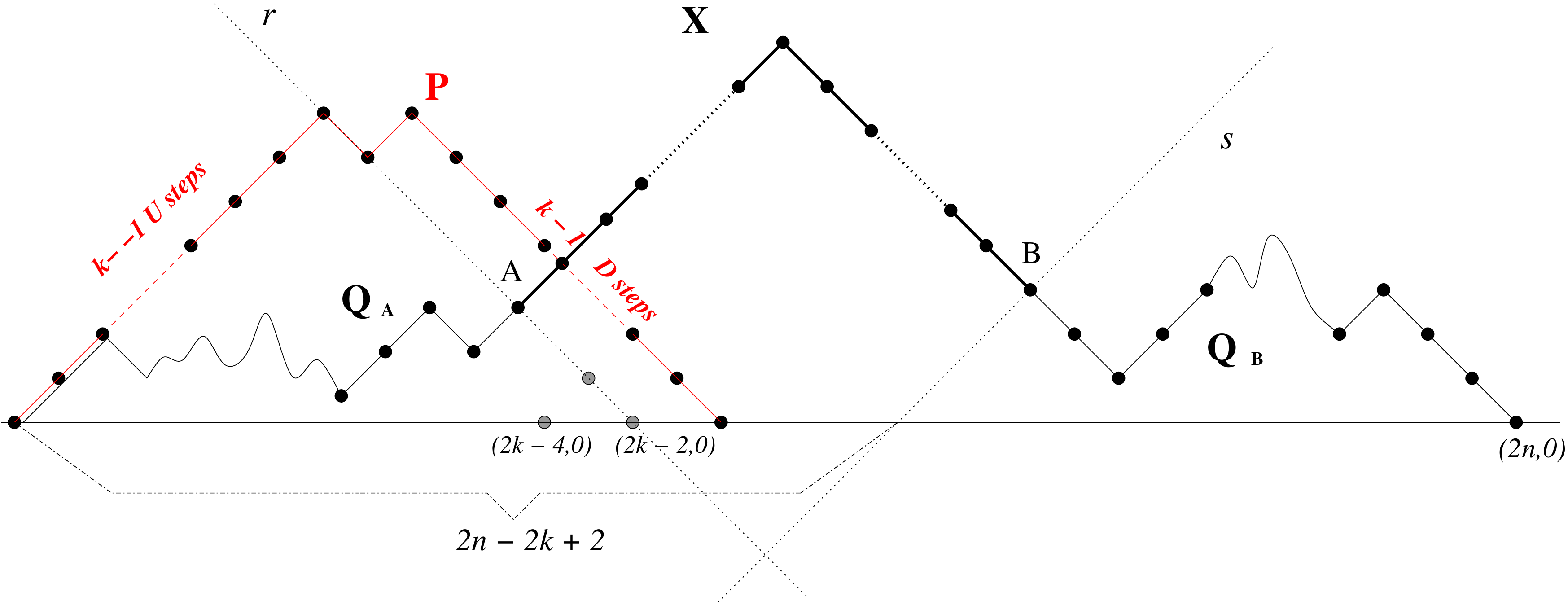}
\end{center}
\vspace{-0.5cm} \caption{Avoiding $U^{k-1}DUD^{k-1}$, with
$n\geq2k-3$} \label{av_UUDUDD_}
\end{figure}

We point out that $Q_A$ has exactly $k-1$ $U$ steps and its last
step is a $U$ step. Analogously, $Q_B$ has exactly $k-1$ $D$ steps
and its first step is a $D$ step. Notice that there is a clear
bijection between the set $\mathcal A$ of Dyck prefixes having
$k-1$ $U$ steps and ending with a $U$ and the set $\mathcal B$ of
Dyck suffixes having $k-1$ $D$ steps and starting with a $D$,
since each element of $\mathcal{B}$ can be read from right to left
thus obtaining an element of $\mathcal{A}$. Moreover,
$\mathcal{A}$ is in bijection with the set of Dyck paths of
semilength $k-1$ (just complete each element of $\mathcal{A}$ with
the correct sequence of $D$ steps), hence $|\mathcal A|=C_{k-1}$.

If we require $Q$ to avoid $P$, then necessarily $X=U^iD^j$, for
suitable $i,j$ (for, if a valley $DU$ occurred in $X$, then $Q$
would contain $P$ since $U^{k-1}$ and $D^{k-1}$ already occur in
$Q_A$ and $Q_B$, respectively). In other words, $A$ and $B$ can be
connected only in one way, using a certain number (possibly zero)
of $U$ steps followed by a certain number (possibly zero) of $D$
steps. Therefore, a path $Q$ avoiding $P$ is essentially
constructed by choosing a prefix $Q_A$ from $\mathcal A$ and a
suffix $Q_B$ from $\mathcal B$, whence:

\begin{equation}\label{ck-1}
d_n(P)=C_{k-1}^2,\quad (\mbox{if}\quad n\geq2k-3).
\end{equation}

\item Suppose now $k+1\leq n<2k-3$ (which means that $r$ and $s$
intersect at height $>1$). Then it can be either $x_A\leq x_B$ or
$x_A>x_B$.

\begin{itemize}

\item[a)] If $x_A\leq x_B$, then we can count all Dyck paths $Q$
avoiding $P$ using an argument analogous to the previous one.
However, in this case the set of allowable prefixes of each such
$Q$ is a proper subset of $\mathcal{A}$. More specifically, we
have to consider only those for which $x_A=k-1,k,k+1,\ldots,n$
(see Figure \ref{av_UUDUDD_2_a}). In other words, an allowable
prefix has $k-1$ $U$ steps and $0,1,2,\ldots$ or $n-k+1$ $D$
steps.
\begin{figure}[h!]
\begin{center}
\includegraphics[scale=.3]{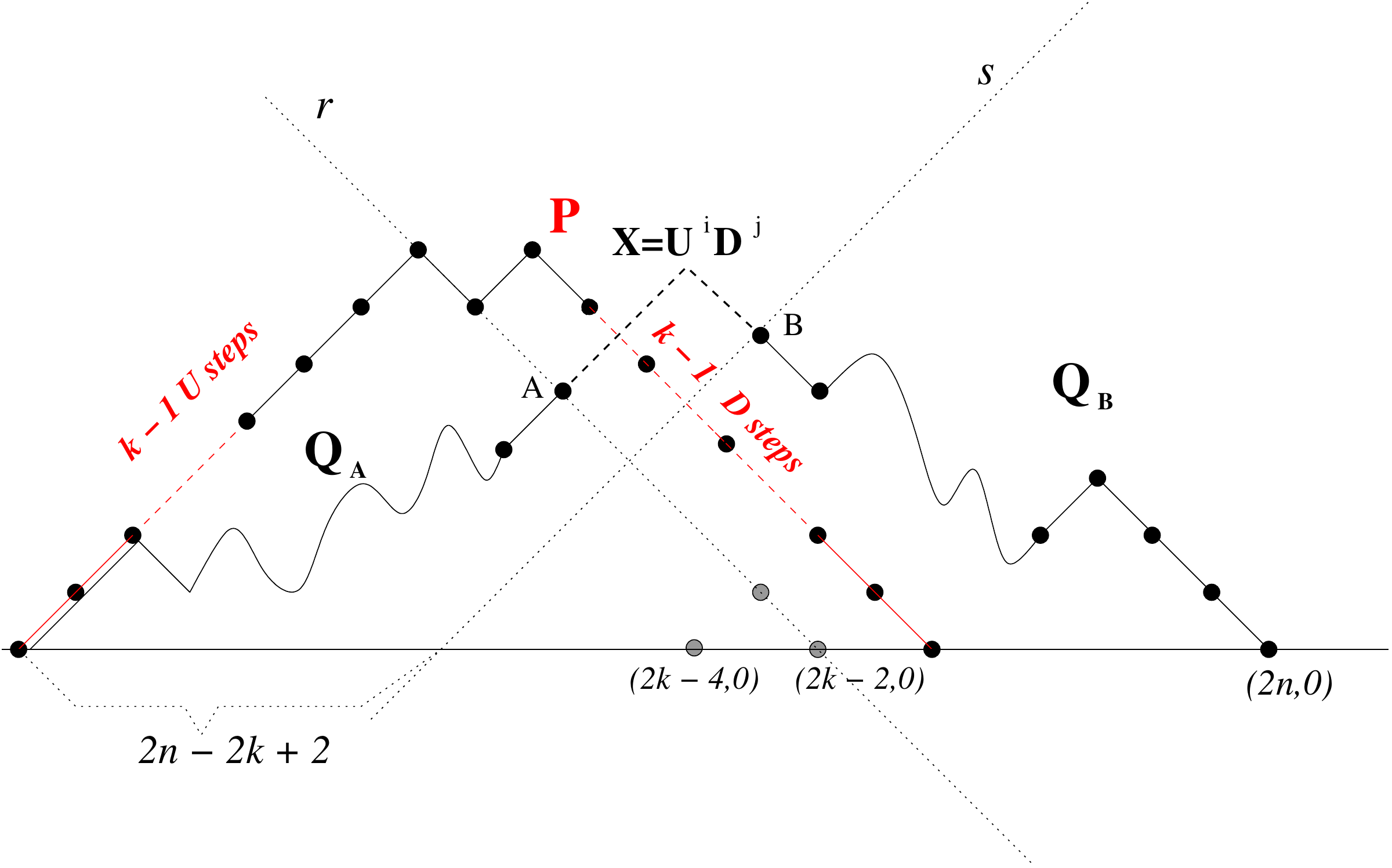}
\end{center}
\vspace{-0.5cm} \caption{Avoiding $U^{k-1}DUD^{k-1}$, with
$x_A\leq x_B$} \label{av_UUDUDD_2_a}
\end{figure}
If $b_{i,j}$ denotes the numbers of Dyck prefixes with $i$ $U$
steps and $j$ $D$ steps ($i\geq j$), then the contribution to $d_n
(P)$ in this case is
$$d_n ^{(1)}(P)=\left(\sum_{j=0}^{n-k+1}b_{k-2,j}\right)^2\quad.$$
The coefficients $b_{i,j}$ are the well-known \emph{ballot
numbers} (sequence A009766 in \cite{Sl}), whose first values are
reported in Table \ref{ballot}.

\item[b)] If $x_A>x_B$, then it is easy to see that $Q$
necessarily avoids $P$, since $A$ clearly occurs after $B$, and so
there are strictly less than $k-1$ $D$ steps from $A$ to $(2n,0)$.
Observe that, in this case, the path $Q$ lies below the profile
drawn by the four lines $y=x$, $r$, $s$ and $y=-x+2n$. In order to
count these paths, referring to Figure \ref{av_UUDUDD_2_b}, just
split each of them into a prefix and a suffix of equal length $n$
and call $C$ the point having abscissa $n$.

\begin{figure}[h!]
\begin{center}
\includegraphics[scale=.3]{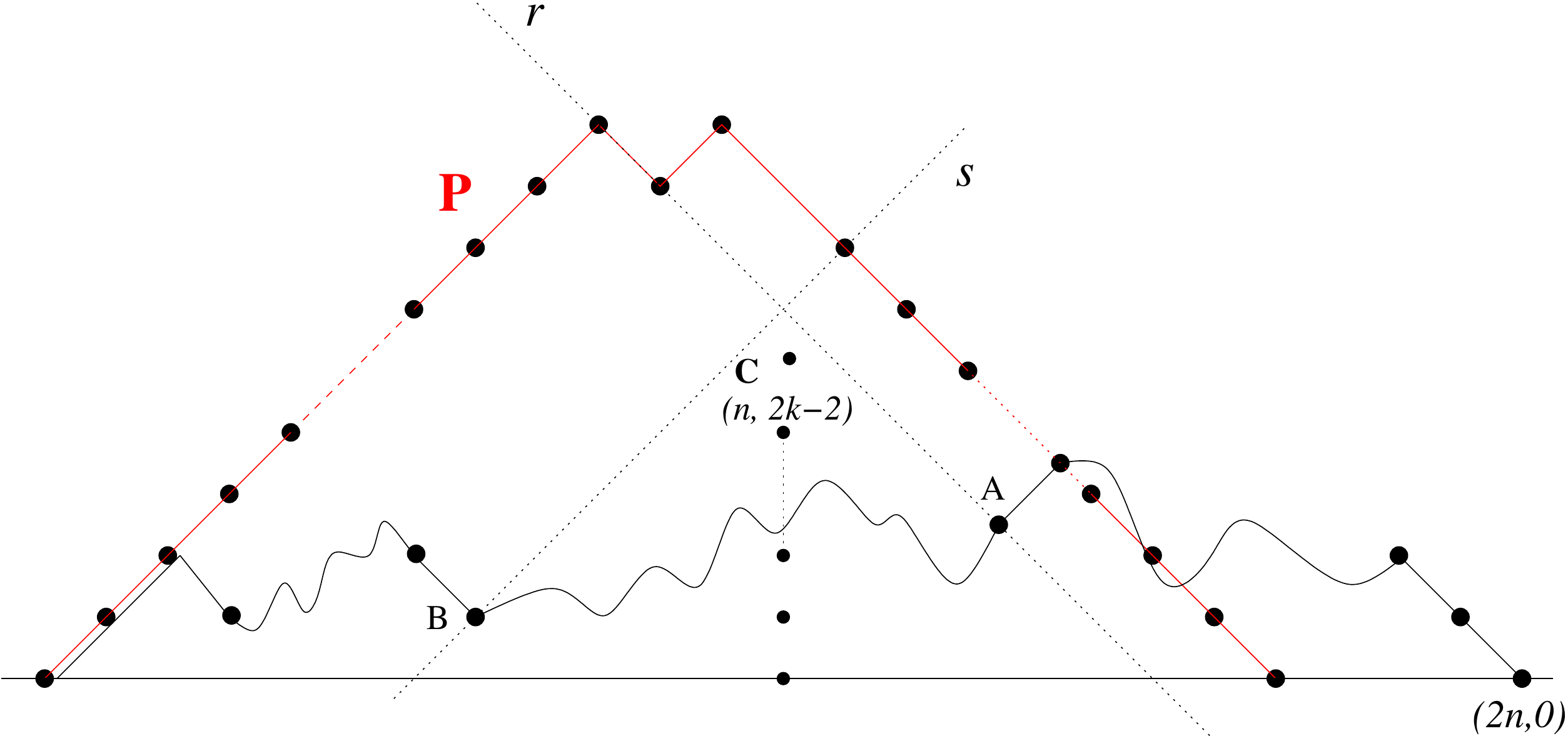}
\end{center}
\caption{Avoiding $U^{k-1}DUD^{k-1}$, with $x_A>x_B$}
\label{av_UUDUDD_2_b}
\end{figure}

Since $C$ must lie under the point where $r$ and $s$ intersect,
then its ordinate $y_C$ equals $-n+2k-2-2t$ with $t\geq 1$ (and
also recalling that $y_C =-n+2k-2-2t\geq 0$). A prefix whose final
point is $C$ has $k-j$ $U$ steps and $n-k+j$ $D$ steps, with
$j\geq 2$. Since, in this case, a path $Q$ avoiding $P$ is
constructed by gluing a prefix and a suffix chosen among
$b_{k-j,n-k+j}$ possibilities $(j\geq 2)$, we deduce that the
contribution to $d_n (P)$ in this case is:
$$d_n ^{(2)}(P)=\sum_{j\geq2}b_{k-j,n-k+j}^2.$$
\end{itemize}

Summing up the two contributions we have obtained in a) and b), we
get:
\begin{align}\label{general_1}
d_n(P)&=d_n ^{(1)}(P)+d_n
^{(2)}(P)\nonumber \\
&=\left(\sum_{j=0}^{n-k+1}b_{k-2,j}\right)^2+\sum_{j\geq2}b_{k-j,n-k+j}^2,\quad
\mbox{if}\quad k+1\leq n<2k-3.
\end{align}

\end{enumerate}

Notice that formula (\ref{general_1}) reduces to the first sum if
$n\geq2k-3$, since in that case $n-k+j>k-j$, for $j\geq2$. We then
have a single formula including both cases 1. and 2.:
\begin{equation}\label{general_2}
d_n(P)=\left(\sum_{j=0}^{n-k+1}b_{k-2,j}\right)^2+\sum_{j\geq2}b_{k-j,n-k+j}^2,
\quad \mbox{if}\quad n\geq k+1\quad.
\end{equation}

Formula (\ref{general_2}) can be further simplified by recalling a
well known recurrence for ballot numbers, namely that
$$
b_{i+1,j}=\sum_{s=0}^{j}b_{i,s}.
$$

Therefore, we get the following interesting expression for $d_n
(P)$ (when $n\geq k+1$) in terms of sums of squares of ballot
numbers along a skew diagonal (see also Tables \ref{ballot} and
\ref{alfa}):

\begin{equation}\label{ultima}
d_n (P)=\sum_{j\geq 1}b_{k-j,n-k+j}^2 .
\end{equation}

\begin{table}
\begin{tabular}{c|cccccccccc}
\backslashbox{i}{j} & 0 & 1 &2&3&4&5&6&7&8&9\\
\hline
0 & 1 &   &   &   &    &    &    &    &&\\
1 & 1 & 1 &   &   &    &    &    &    &&\\
2 & 1 & 2 &  2&   &    &    &    &    &&\\
3 & 1 & 3 &  5&  5&    &    &    &    &&\\
4 & 1 & 4 &  9& 14&  14&    &    &    &&\\
5 & 1 & 5 & 14& 28&  \textbf{42}&  42&    &&    \\
6&\multicolumn{1}{>{\columncolor[gray]{.7}}c}{\emph{1}} &
\multicolumn{1}{>{\columncolor[gray]{.7}}c}{\emph{6}} &
\multicolumn{1}{>{\columncolor[gray]{.7}}c}{\emph{20}}&
\textbf{48}&  90& 132& 132&&&    \\
7 & 1 & 7 & \textbf{27}& 75& 165& 297& 429& 429&&\\
8 & 1 & 8 & 35 & 110 & 275 & 572 & 1001 & 1430 & 1430\\
9 & 1 & 9 & 44 & 154 & 429 & 1001 & 2002 & 3432 & 4862 & 4862\\
\end{tabular}
\caption{The sum of the gray entries gives the bold entry in the
line below. The sum of the squares of the bold entries gives an
appropriate element of Table \ref{alfa}.} \label{ballot}
\end{table}

\begin{table}
{\footnotesize
\begin{tabular}{c|cccccccccccccccc}

\backslashbox{$k$}{$n$} & 0 & 1 & 2 & 3 & 4 & 5 & 6 & 7 &8&9&10&11&12&13&\ldots\\
\hline
1 & 0 & 0 & 0 & 0 & 0 & 0 & 0 & 0 & 0& 0& 0& 0& 0& 0& \ldots \\
2 & 1 & 1 & 1 & 1 & 1 & 1 & 1 & 1 & 1& 1& 1& 1& 1& 1& \ldots \\
3 & 1 & 1 & 2 &  {\slshape 4} & 4 & 4 & 4 & 4 & 4& 4& 4& 4& 4& 4& \ldots \\
4 & 1 & 1 & 2 & 5 & {\slshape 13} & 25 & 25 & 25 & 25& 25& 25& 25& 25& 25 & \ldots \\
5 & 1 & 1 & 2 & 5 & 14 & {\slshape 41} & \textbf{106} & 196 & 196 & 196 & 196 & 196 & 196 & 196 & \ldots    \\

6 & 1 & 1 & 2 & 5 & 14 & 42 & {\slshape 131} & \textbf{392} & \textbf{980} & 1764 & 1764 & 1764 & 1764 & 1764 & \ldots \\

7 & 1 & 1 & 2 & 5 & 14 & 42 & 132 & {\slshape 428} & \textbf{1380} & \textbf{4068} & \textbf{9864} & 17424 & 17424 & 17424 & \ldots   \\

8 & 1 & 1 & 2 & 5 & 14 & 42 & 132 & 429 & {\slshape 1429} & \textbf{4797} & \textbf{15489} & \textbf{44649} & \textbf{105633} & 184041 & \ldots\\

9 & 1 & 1 & 2 & 5 & 14 & 42 & 132 & 429 & 1430 & {\slshape 4861} & \textbf {16714} & \textbf{56749} & \textbf{181258} & \textbf{511225} & \ldots\\
\end{tabular}
} \vspace{1cm} \caption{Number of Dyck paths of semilength $n$
avoiding $U^{k-1}DUD^{k-1}$. Entries in boldface are the
nontrivial ones ($k+1\leq n<2k-3$).} \label{alfa}
\end{table}

Therefore we obtain in particular:
$$d_n (UUDUDD)=4,\textnormal{ when $n\geq 3$.}$$

\subsection{The pattern $U^k D^k$}

The case $P=U^kD^k$ is very similar to the previous one. We just
observe that, when $x_A\leq x_B$, the two points $A$ and $B$ can
be connected only using a sequence of $D$ steps followed by a
sequence of $U$ steps. This is possible only if $n\leq 2k-2$,
which means that $r$ and $s$ do not intersect below the $x$-axis.
Instead, if $n\geq 2k-1$, $Q$ cannot avoid $P$. Therefore we get
(see also Table \ref{beta}):

$$
d_n(P)= \left\{
\begin{array}{cc}
0 & \mbox{if $n\geq 2k-1$};\\
\sum_{j\geq1}b_{k-j,n-k+j}^2 & \mbox{otherwise}.
\end{array}
\right.
$$

In particular, we then find:

\begin{itemize}

\item[-] $d_n (UUDD)=0$, when $n\geq 3$;

\item[-] $d_n (UUUDDD)=0$, when $n\geq 5$.

\end{itemize}

\begin{table}
{\footnotesize
\begin{tabular}{c|cccccccccccccccc}
\backslashbox{$k$}{$n$} &
    0 & 1 & 2 & 3 & 4 & 5 & 6 & 7 & 8 & 9 & 10 & 11 & 12 & 13 & \ldots\\
\hline
1 & 0 & 0 & 0 & 0 & 0 & 0 & 0 & 0 & 0 & 0 & 0 & 0 & 0 & 0 & \ldots \\

2 & 1 & 1 & 1 & 0 & 0 & 0 & 0 & 0 & 0 & 0 & 0 & 0 & 0 & 0 & \ldots \\

3 & 1 & 1 & 2 & \slshape{4} & 4 & 0 & 0 & 0 & 0 & 0 & 0 & 0 & 0 & 0 & \ldots \\

4 & 1 & 1 & 2 & 5 & \slshape{13} & 25 & 25 & 0 & 0 & 0 & 0 & 0 & 0 & 0 & \ldots \\

5 & 1 & 1 & 2 & 5 & 14 & \slshape{41} & \textbf{106} & 196 & 196 & 0 & 0 & 0 & 0 & 0 & \ldots    \\

6 & 1 & 1 & 2 & 5 & 14 & 42 & \slshape{131} & \textbf{392} & \textbf{980} & 1764 & 1764 & 0 & 0  & 0 & \ldots  \\

7 & 1 & 1 & 2 & 5 & 14 & 42 & 132 & \slshape{428} & \textbf{1380} & \textbf{4068} & \textbf{9864} & 17424 & 17424 & 0 &\ldots  \\

8 & 1 & 1 & 2 & 5 & 14 & 42 & 132 & 429 & \slshape{1429} & \textbf{4797} & \textbf{15489} & \textbf{44649} & \textbf{105633} & 184041 & \ldots\\

9 & 1 & 1 & 2 & 5 & 14 & 42 & 132 & 429 & 1430 & \slshape{4861} & \textbf{16714} & \textbf{56749} & \textbf{181258} & \textbf{511225} & \ldots\\
\end{tabular}
} \vspace{1cm} \caption{Number of Dyck paths of semilength $n$
avoiding $U^{k}D^{k}$. Entries in boldface are the nontrivial ones
($k+1\leq n<2k-3$).} \label{beta}
\end{table}

\subsection{The pattern $U^{k-1}D^{k-1}UD$}

This is by far the most challenging case.

Let $Q$ be a Dyck path of length $2n$ and $P=U^{k-1}D^{k-1}UD$. If
$Q$ avoids $P$, then there are two distinct options: either $Q$
avoids $U^{k-1}D^{k-1}$ or $Q$ contains such a pattern. In the
first case, we already know that $d_n (U^{k-1}D^{k-1})$ is
eventually equal to zero. So, for the sake of simplicity, we will
just find a formula for $d_n (P)$ when $n$ is sufficiently large,
i.e. $n\geq 2k-3$. Therefore, for the rest of this section, we
will suppose that $Q$ contains $U^{k-1}D^{k-1}$.

\bigskip

\begin{figure}[h!]
\begin{center}
\includegraphics[scale=.3]{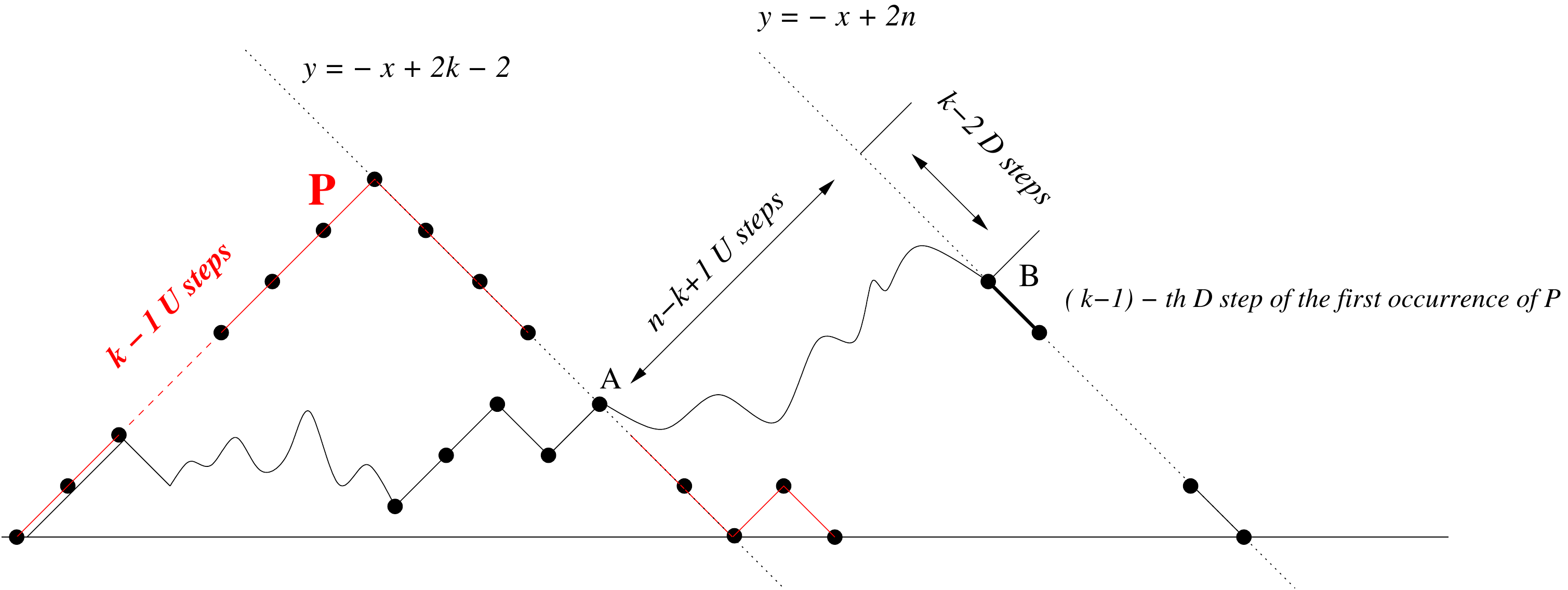}
\end{center}
\caption{A path $Q$ avoiding $P=U^{k-1}D^{k-1}UD$}
\label{challenging_1}

\end{figure}

The $(k-1)$-th $D$ step of the first occurrence of
$U^{k-1}D^{k-1}$ in $Q$ lies on the line having equation
$y=-x+2n$. This is due to the fact that $Q$ has length $2n$ and
there cannot be any occurrence of $UD$ after the first occurrence
of $U^{k-1}D^{k-1}$. The path $Q$ touches the line of equation
$y=-x+2k-2$ for the first time with the end point $A$ of its
$(k-1)$-th $U$ step. After that, the path $Q$ must reach the
starting point $B$ of the $(k-1)$-th $D$ step occurring after $A$.
Finally, a sequence of consecutive $D$ steps terminates $Q$ (see
Figure \ref{challenging_1}). Therefore, $Q$ can be split into
three parts: the first part, from the beginning to $A$, is a Dyck
prefix having $k-1$ $U$ steps and ending with a $U$ step; the
second part, from $A$ to $B$, is a path using $n-k+1$ $U$ steps
and $k-2$ $D$ steps; and the third part, from $B$ to the end, is a
sequence of $D$ steps (whose length depends on the coordinates of
$A$). However, both the first and the second part of $Q$ have to
obey some additional constraints.

The height of the point $A$ (where the first part of $Q$ ends)
must allow $Q$ to have at least $k-1$ $D$ steps after $A$. Thus,
the height of $A$ plus the number of $U$ steps from $A$ to $B$
minus the number of $D$ steps from $A$ to $B$ must be greater than
or equal to 1 (to ensure that the pattern $U^{k-1}D^{k-1}$ occurs
in $Q$). Hence, denoting with $x$ the maximum number of $D$ steps
which can occur before $A$, either $x=k-2$ or the following
equality must be satisfied:
$$
(k-1)-x+(n-k+1)-(k-2)=1.
$$
Therefore, $x=\min \{ n-k+1,k-2\}$. Observe however that, since we
are supposing that $n\geq 2k-3$, we always have $x=k-2$.

\begin{figure}[h!]
\begin{center}
\includegraphics[scale=.3]{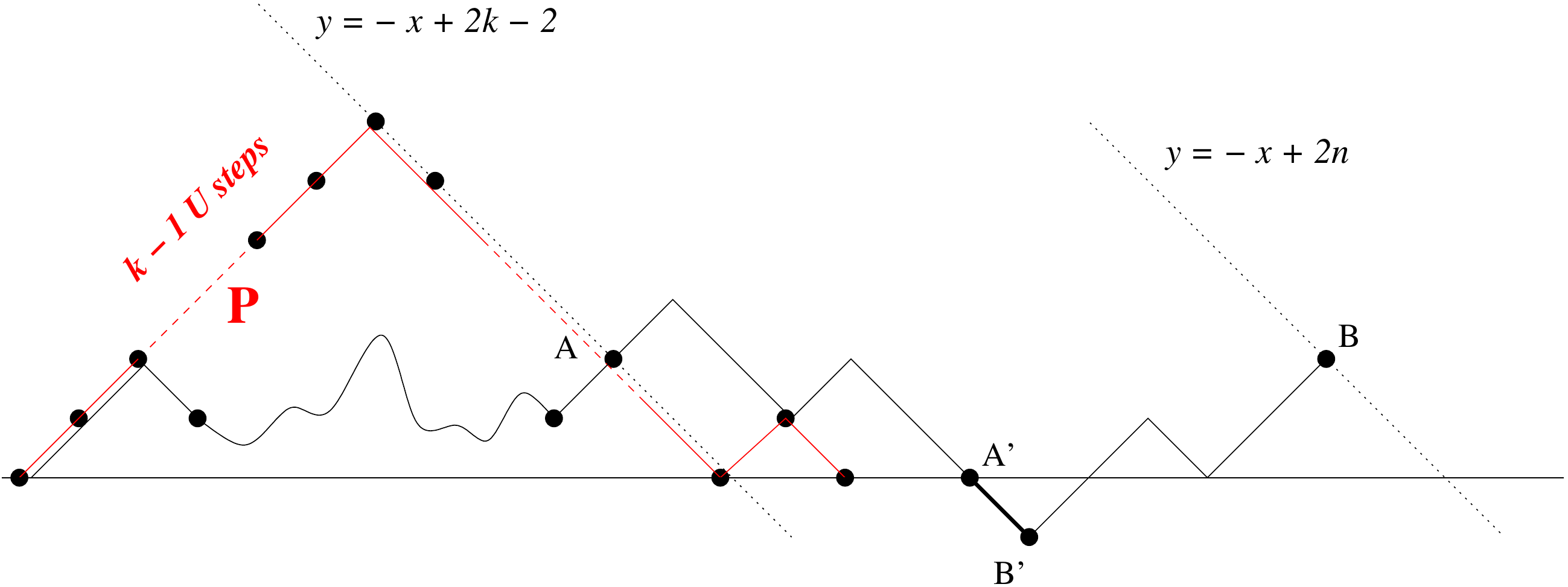}
\end{center}
\caption{A forbidden subpath from $A$ to $B$.}
\label{challenging_2}
\end{figure}

Concerning the part of $Q$ between $A$ and $B$, since we have to
use $n-k+1$ $U$ steps and $k-2$ $D$ steps, there are ${n-1\choose
k-2}$ distinct paths connecting $A$ and $B$. However, some of them
must be discarded, since they fall below the $x$-axis. In order to
count these ``bad" paths, we split each of them into two parts.
Namely, if $A'$ and $B'$ are the starting and ending points of the
first (necessarily $D$) step below the $x$-axis, the part going
from $A$ to $A'$, and the remaining part (see Fig.
\ref{challenging_2}).
It is not too hard to realize that the number of possibilities we
have to choose the first part is given by a ballot number
(essentially because, reading the path from right to left, we have
to choose a Dyck prefix from $A'$ to $A$), whereas the number of
possibilities we have to choose the second part is given by a
binomial coefficient (essentially because, after having discarded
the step starting at $A'$, we have to choose an unrestricted path
from $B'$ to $B$). After a careful inspection, we thus get to the
following expression for the total number $d_n (P)$ of Dyck paths
of semilength $n\geq 2k-3$ avoiding $P$:
\begin{align}\label{bad}
d_n (P)=&{n-1\choose k-2}C_{k-1}\nonumber \\
&-\sum_{s=2}^{k-2}b_{k-2,s}\cdot \left(
\sum_{i=0}^{s-2}b_{k-3-i,s-2-i}{n-k-s+3+2i\choose i}\right) .
\end{align}

Formula (\ref{bad}) specializes to the following expressions for
low values of $k$ (see also Table \ref{difficile}):

\begin{itemize}

\item[-] when $k=3$, $d_n (P)=2n-2$ for $n\geq 3$;

\item[-] when $k=4$, $d_n (P)=\frac{5n^2 -15n+6}{2}$ for $n\geq
5$;

\item[-] when $k=5$, $d_n (P)=\frac{14n^3 -84n^2 +124n-84}{6}$ for
$n\geq 7$.

\end{itemize}

\begin{table}
\begin{tabular}{c|cccccccccc}
\backslashbox{$k$}{$n$} &
    0 & 1 & 2 & 3 & 4 & 5 & 6 & 7 & 8 & 9 \\
\hline
1 & 1 & 1 & 0 & 0 & 0 & 0 & 0 & 0 & 0 & 0 \\

2 & 1 & 1 & 1 & 1 & 1 & 1 & 1 & 1 & 1 & 1 \\

3 & 1 & 1 & 2 & 4 & 6 & 8 & 10 & 12 & 14 & 16 \\

4 & 1 & 1 & 2 & 5 & 13 & 28 & 48 & 73 & 103 & 138 \\

5 & 1 & 1 & 2 & 5 & 14 & 41 & 110 & 245 & 450 & 739 \\

6 & 1 & 1 & 2 & 5 & 14 & 42 & 131 & 397 & 1069 & 2427 \\
\end{tabular}
\caption{Avoiding $U^{k-1}D^{k-1}UD$} \label{difficile}
\end{table}

\section{Some remarks on the asymptotics of pattern avoiding Dyck paths}

In this final section we collect some thoughts concerning the
asymptotic behavior of integer sequences counting pattern-avoiding
Dyck paths. Unlike the case of permutations, for Dyck paths it
seems plausible that a sort of ``master theorem'' exists, at least
in the case of single avoidance. This means that all the sequences
which count Dyck paths avoiding a single pattern $P$ have the same
asymptotic behavior (with some parameters, such as the leading
coefficient, depending on the specific path $P$). We have some
computational evidence which leads us to formulate a conjecture,
whose proof we have not been able to complete, and so we leave it
as an open problem.

\bigskip

Let $P$ denote a fixed Dyck path of semilength $x$. We are
interested in the behavior of $d_n (P)$ when $n\rightarrow
\infty$. Our conjecture is the following:

\bigskip

\textbf{Conjecture.}\quad \emph{Suppose that $P$ starts with $a$ U
steps and ends with $b$ D steps. Then, setting $k=2x-2-a-b$, we have
that $d_n (P)$ is asymptotic to
$$\frac{\alpha_P \cdot C_a \cdot
C_b}{k!}n^{k},$$ where $C_m$ denotes the $m$-th Catalan numbers
and $\alpha_P$ is the number of saturated chains in the Dyck
lattice of order $x$ (see \cite{FP}) from $P$ to the maximum $U^x
D^x$.}

\bigskip

\begin{figure}[h!]
\begin{center}
\includegraphics[scale=.3]{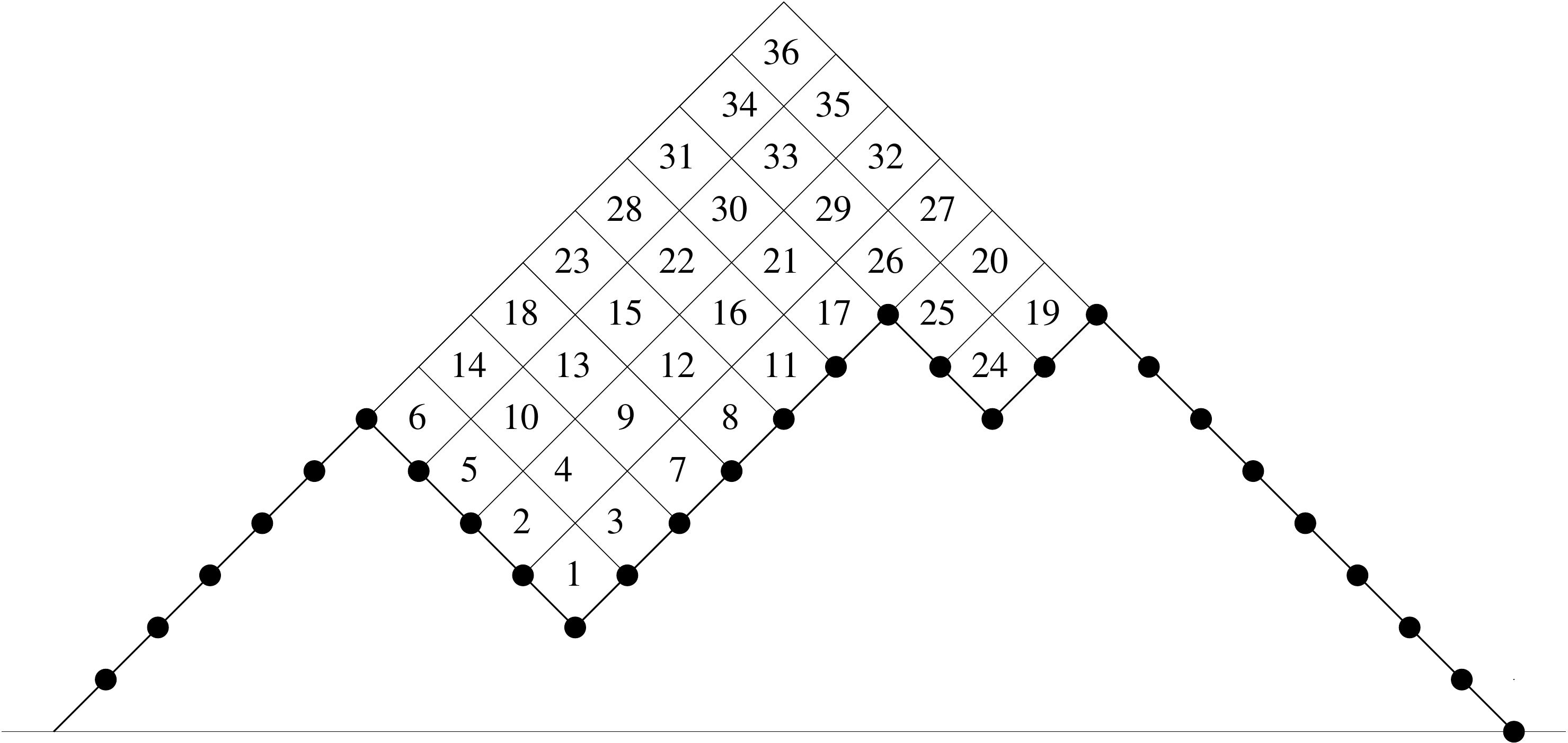}
\end{center}
\caption{An instance of a standard Young tableau determined by a
Dyck path.} \label{young}

\end{figure}

Equivalently, $\alpha_P$ is the number of standard Young tableaux
whose Ferrers shape is determined by the region delimited by the
path $P$ and the path $U^x D^x$, as shown in Figure \ref{young}.

In the above conjecture, the only parts of the formula we are able
to justify are the coefficients $C_a$ and $C_b$. Indeed, suppose
that $Q$ is a Dyck path of semilength $n$, with $n$ very large.
Then we can consider the minimum prefix $Q_{pref}$ of $Q$
containing exactly $a$ $U$ steps and the minimum suffix $Q_{suff}$
of $Q$ containing exactly $b$ $D$ steps. They certainly exist, due
to the hypothesis that $n$ is very large. As we have already shown
in the previous section, the number of Dyck prefixes having $a$
$U$ steps and ending with $U$ is precisely equal to $C_a$. Of
course, an analogous fact holds for suffixes as well.

\bigskip

We close our paper with some further conjectures concerning the
order structure of the Dyck pattern poset.

\begin{itemize}

\item What is the M\"obius function of the Dyck pattern poset
(from the bottom element to a given path? Of a generic interval?)?

\item How many (saturated) chains are there up to a given path? Or
in a general interval?

\item Does there exist an infinite antichain in the Dyck pattern
poset?

\end{itemize}

The last conjecture has been suggested by an analogous one for the
permutation pattern poset which has been solved in the affirmative
(see \cite{SB} and the accompanying comment). In the present
context we have no intuition on what could be the answer, though
we are a little bit less optimistic than in the permutation case.


\begin{thebibliography}{99}

\bibitem[Bj]{Bj} A. Bj\"orner,\quad \emph{The M\"obius function of subword
order},\quad Invariant theory and tableaux (Minneapolis, MN,
1988), 118–124, IMA Vol. Math. Appl., 19, Springer, New York,
1990.

\bibitem[Bu]{Bu} A. Burstein,\quad \emph{Enumeration of words with forbidden
patterns},\quad PhD thesis, University of Pennsylvania, 1998.

\bibitem[DPTW]{DPTW} M. Dairyko, L. Pudwell, S. Tyner, C.
Wynn,\quad \emph{Non-contiguous pattern avoidance in binary
trees},\quad Electron. J. Combin.,\quad 19(3) (2012) \#P22.

\bibitem[FP]{FP} L. Ferrari, R. Pinzani,\quad \emph{Lattices of lattice
paths},\quad J. Statist. Plann. Inference,\quad 135 (2005) 77--92.

\bibitem[Gi]{Gi} Samuele Giraudo,\quad \emph{Balanced binary trees in the Tamari
lattice},\quad Discrete Math. Theoret. Comput. Sci. Proc.,\quad AN
(2010) 593--604.

\bibitem[Go]{Go} A. Goyt,\quad \emph{Avoidance of partitions of a three element
set},\quad Adv. Appl. Math.,\quad 41 (2008) 95--114.

\bibitem[Kl]{Kl} M. Klazar,\quad \emph{On $abab$-free and $abba$-free sets
partitions},\quad European J. Combin.,\quad 17 (1996) 53--68.

\bibitem[Kn]{Kn} D. Knuth,\quad \emph{The art of computer programming}, Vol. 1,
Addison Wesley, Boston, 1968.

\bibitem[R]{R} E. Rowland,\quad \emph{Pattern avoidance in binary
trees},\quad J. Combin. Theory Ser. A,\quad 117 (2010) 741--758.

\bibitem[Sa]{Sa} B. E. Sagan,\quad \emph{Pattern avoidance in set
partitions},\quad Ars Combin.,\quad 94 (2010) 79--96.

\bibitem[Sl]{Sl} N. J. A. Sloane,\quad \emph{The On-Line Encyclopedia of Integer Sequences},\quad
electronically available at {\tt oeis.org}.

\bibitem[SB]{SB} D. A. Spielman, M. B\'{o}na,\quad \emph{An infinite antichain of
permutations},\quad Electron. J. Combin.,\quad 7 (2000) \#N2 (4
pp.).

\end{thebibliography}
\end{document}